\documentclass[12pt]{article}
\input{epsf.sty}
\usepackage{hyperref}
\usepackage{graphicx}
\usepackage{amsmath, amssymb}
\usepackage[arrow, matrix, curve]{xy}
\newtheorem {thm}{Theorem}[section]
\newtheorem {prop}[thm]{Proposition}

\newtheorem {defn}[thm]{Definition}

\newcommand{\qed}{\nobreak \ifvmode \relax \else
      \ifdim\lastskip<1.5em \hskip-\lastskip
      \hskip1.5em plus0em minus0.5em \fi \nobreak
      \vrule height0.75em width0.5em depth0.25em\fi}

\def\Cox{\hfill \Box}

\def\N{{\Bbb N}}

\def\Z{{\Bbb Z}}
\def\R{{\Bbb R}}

\def\e{{\varepsilon}}

\def\D{\Delta}
\def\a{\alpha}

\def\b{\beta}

\def\d{\delta}

\def\e{\varepsilon}

\def\phi{\varphi}

\def\g{\gamma}

\def\l{\lambda}

\def\r{\rho}

\def\s{\sigma}

\def\D{\Delta}

\def\L{\Lambda}

\def\G{\Gamma}

\def\O{{\Omega}}

\def\T{\T}

\def\CC{{\cal C}}

\def\PP{{\cal P}}

\begin{document}
\title{Sharp thresholds for Gibbs-non-Gibbs transition in the fuzzy Potts model with a Kac-type interaction}

\author{
Benedikt Jahnel
\footnote{ Weierstrass Institute Berlin, Mohrenstr. 39, 10117 Berlin, Germany,
\newline
 \texttt{Benedikt.Jahnel@wias-berlin.de}, 
\newline
\texttt{http://www.wias-berlin.de/people/jahnel/ }}
 \, and  Christof K\"ulske
\footnote{ Ruhr-Universit\"at   Bochum, Fakult\"at f\"ur Mathematik, D44801 Bochum, Germany,
\newline
\texttt{Christof.Kuelske@ruhr-uni-bochum.de}, 
\newline
\texttt{http://www.ruhr-uni-bochum.de/ffm/Lehrsttuehle/Kuelske/kuelske.html
/$\sim$kuelske/ }}\, 
\,  
}

\maketitle

\begin{abstract}
We investigate the Gibbs properties of the fuzzy Potts model on the $d$-dimensional 
torus with Kac interaction. 
We use a variational approach for profiles inspired by that of Fern\'andez, 
den Hollander and Mart\'inez \cite{FeHoMa14} for their study of the Gibbs-non-Gibbs transitions of a 
dynamical Kac-Ising model on the torus.  
As our main result, we show that the mean-field thresholds dividing Gibbsian from non-Gibbsian behavior are sharp in the fuzzy Kac-Potts model with class size unequal two.  
On the way to this result we prove a large deviation principle  
for color profiles with diluted total mass densities and use monotocity arguments. 

\end{abstract}




\vfill\eject

\section{Introduction}

In previous years we have seen a number of measures describing systems with interacting components appearing in mathematical statistical mechanics which have lost the Gibbs property as a result of a transformation \cite{En12,Ny08,EnFeHoRe02,EnFeHoRe10,RoMaNeSc13}. 
Such a loss is indicated by the failure of continuity of conditional probabilities at a given site, when the conditioning is varied away from this site. 
Interesting sources of non-Gibbsian behavior include time evolutions or deterministic transformations which reduce the complexity of the local state space. 
A prototypical example of a system of the second type is the fuzzy Potts model (fuzzy PM) \cite{HaKu04,Po52,KuRo14,JaKuRuWe14,Ha03,Ba10}. 
It is obtained from the ordinary PM by partitioning the local state space $\{1,2,\dots,q\}$ into subclasses and observing the Potts distribution after identification of the spin-values inside the subclasses. 

It has been noted in some cases for mean-field models \cite{Ku03,JaKuRuWe14,HoReZu15} when the appropriate notion of mean-field Gibbsianness is employed, the question of continuity can be reduced to variational problems. For systems for which lattice results and mean-field results are available it turns out that these results are often in a striking parallel \cite{KuNy07,EnFeHoRe10}. It is an open challenge to understand this relation better. 

One way to approach the relation between the lattice and mean-field is via Kac models (KM) \cite{BeMoOrSaTr12,Co87,Co89,EiEl83,MaOrPrTr94,BoKu05} in which there is a parameter which makes the interaction long-range but a spatial structure remains. 

The first rigorous result relating Gibbs properties of a KM to that of a mean-field model was obtained in \cite{FeHoMa14} in the case of independent time evolutions from an initial Kac-Ising model. The relation between a spatial model and a mean-field model was set up as follows. 
The authors put the model on a torus in $d$ dimensions, with spins sitting on a grid of spacing $1/n$,  
and looked at a single-site conditional probability in the large $n$-limit. 
The limiting object they studied then was a specification kernel giving the dependence of a single-site probability 
as a function of a magnetization profile. 
The existence of the limiting kernel and properties of its approach along volume sequences were established using a combination of a large deviation principle (LDP) in equilibrium for the Ising model \cite{Co87}, a path LDP, and techniques from hydrodynamic limits. It was not possible to give sharp parameter values for the Gibbs-non-Gibbs (GnG) transition but sufficient conditions on time and initial temperature values to be non-Gibbsian could be provided.

\medskip
In our present study of the fuzzy Kac-Potts model (fuzzy KPM) we ask related questions. Our main result is Theorem \ref{HaKu-Kac} where we provide 
precise threshold values dividing Gibbsian and non-Gibbsian behavior. To our knowledge this is the first sharp 
result for GnG in a KM.

\subsection{Strategy of proof and further results}
The Hamiltonian of the KPM can be written in terms of an 
empirical color distribution field and we start by noting 
a LDP for the empirical color distribution field as the grid on the torus shrinks. 
The minimizers of the rate function for this LDP provide us with the equilibrium phases, and it is easy to see that the absolute minimizers must be flat (spatially homogeneous). Therefore the critical value for phase transitions in the KPM is given by the corresponding mean-field result (the Ellis-Wang Theorem \cite{ElWa89}). 

Next, to investigate the Gibbsian properties of the \textit{fuzzy} model we analyse limiting expressions for the single-site conditional probabilities (the specification kernel). 
The idea to prove equality of critical parameters dividing GnG in mean-field with the corresponding critical parameters in the KPM is then to make rigorous the statement that there are no worse conditionings than spatially homogeneous conditionings for fuzzy classes of size unequal two. 
As an intermediate step we prove a LDP for color profiles for a spatially diluted KPM in Proposition \ref{DiLDP}. This and the corresponding non-homogeneous variational problems are interesting in their own right. 
We relate the specification kernel to solutions of such variational problems 
where the dilutions are prescribed by the conditioning profile. Finally this is supplemented by monotonicity arguments in the dilution to show sharpness of the mean-field values for the KM.
%
%
%
%
%
%
%

\subsection{Acknowledgment}
This work is supported by the Sonderforschungsbereich SFB $|$ TR12-Symmetries and Universality in Mesoscopic Systems. Christof K\"ulske thanks F.~Comets, R.~Fern\'andez, F.~den Hollander and J.~Mart\'{\i}nez  for stimulating discussions.

\section{Model and main results}
\subsection{The Kac-Potts model}

Let $\mathbb{T}^d :=\R^d/\Z^d$ be the $d$-dimensional unit torus. For $n\in \N$, let $\mathbb{T}^d_n$ be the $(1/n)$-discretization of $\mathbb{T}^d$ defined by $\mathbb{T}^d_n:=\D^d_n/n$, with $\D^d_n:=\Z^d /n\Z^d$ the discrete torus of size $n$. For $n\in\N$, let $\O_n:=\{1,\dots,q\}^{\D^d_n}$ be the set of \textit{Potts-spin configurations} on $\D^d_n$. We will call elements of $\{1,\dots,q\}$ \textit{colors}. The energy of the configuration $\s := (\s(x))_{x\in\D^d_n}\in\O_n$ is given by the \textit{Kac-type Hamiltonian} 
\begin{equation}
H_n(\s):=-\frac{1}{n^d}\sum_{x,y\in\D^d_n}J(\frac{x-y}{n})1_{\s(x)=\s(y)}, \hspace{1cm}\s\in\O_n
\end{equation}
where $0\leq J\in C(\mathbb{T}^d)$ is a continuous \textit{interaction-functions} on $\mathbb{T}^d$ which is symmetric and $\int dv J(v)=1$. The Gibbs measure associated with $H_n$ is given by
\begin{equation}\label{KPM}
\mu_n(\s):=\frac{1}{Z_n}\exp(-\b H_n(\s)),\hspace{1cm}\s\in\O_n
\end{equation}
with $\b\in[0,\infty)$ the inverse temperature and $Z_n$ the normalizing partition sum.

\bigskip
We are interested in the large $n$-limit for $\mu_n$ and prepare the analysis by rewriting the Hamiltonian in terms of density profiles. More precisely, for $\L\subset\D^d_n$ let $\pi_{\L}: \O_n\mapsto\PP(\mathbb{T}^d_n\times\{1,\dots,q\})\subset\PP(\mathbb{T}^d\times\{1,\dots,q\})$ be the \textit{empirical color measure vector} or \textit{color profiles }of $\s$ inside the volume $\L$ defined by
\begin{equation*}
\pi_{\L}^\s:=\frac{1}{|\L|}\Bigl(\sum_{x\in\L}1_{\s(x)=1}\d_{x/n},\dots,\sum_{x\in\L}1_{\s(x)=q}\d_{x/n}\Bigr)^{T}
\end{equation*} 
where $\d_u$ is the point measure at $u\in\mathbb{T}^d$. In the sequel we use notation $\PP_n:=\PP(\mathbb{T}^d_n\times\{1,\dots,q\})$ and $\PP:=\PP(\mathbb{T}^d\times\{1,\dots,q\})$.
For any $\nu\in \PP$ we will write $\nu[a]$ to indicate the evaluation of $\nu$ at a color $a\in\{1,\dots,q\}$, in other words, $\nu[a]$ is the spatial profile of sites with color $a$. 
In particular, for $x\in \L$, $\pi_{\L}^\s[a](x/n)=|\L|^{-1}1_{\s(x)=a}$.
%
%
%
%
%

Let $u\in\mathbb{T}^d$, then for the \textit{color profile perforated at} $u\in\mathbb{T}^d$ we write $\pi_n^{(u)}:=\pi_{\D_n^d\setminus\lfloor nu\rfloor}$ where $\lfloor nu\rfloor$ denotes the lower-integer part of $nu$. 
Further we abbreviate 
$\mathcal{M}_n:=\pi_n(\O_n)\subset\PP_n$ and $\mathcal{M}_n^u:=\pi_n^{(u)}(\O_n)\subset\PP$ for the sets of possible profiles of mesh-size $n$ and possible profiles of mesh-size $n$ perforated at site $u$. 

We equip $\PP$ and the indicated subspaces with the weak topology, i.e.~the topology corresponding to convergence of continuous functions \linebreak $f\in C(\mathbb{T}^d\times\{1,\dots,q\},\R)=:\CC$. 
This convergence can be metrized in the usual way (see for example \cite[page 235]{Ba78}) by choosing a dense set of functions $( f_j)_{j\in\N}\subset \CC$ and setting 
\begin{equation}\label{weakmetric}
\begin{split}
d(\mu,\nu):=\sum_{j=1}^\infty 2^{-j}\frac{|\mu( f_j)-\nu( f_j)|}{1+|\mu( f_j)-\nu( f_j)|}.
\end{split}
\end{equation}
Moreover since $\mathbb{T}^d\times\{1,\dots,q\}$ is compact and Polish also $(\PP,d)$ is compact and Polish.
Notice that $\s\in\O_n$ determines $\pi_{n}^\s\in\PP_n$ and vice versa.
%
%

\medskip
Using color profiles, we can rewrite the Hamiltonian as
\begin{equation}
H_n(\s)=-n^d \sum_{a=1}^qF(\pi_{n}^\s[a])
\end{equation}
with $F(\nu[a]):=\langle J\ast \nu[a],\nu[a]\rangle=\int\int \nu[a](du)\nu[a](dv)J(u-v)$.
We will be interested in weak limits of color profiles in $\PP$, especially those  having $q$-dimensional Lebesgue densities of the form $\nu=\a\l=(\a[1]\l,\dots,\a[q]\l)^T$ with $\a\in B$ where 
\begin{equation}
\begin{split}
B:=\{\a=(\a[1],\dots,\a[q])^T: \text{ }&0\leq\a[a]\in L^\infty(\mathbb{T}^d,\l)\cr
&\text{ with }
\sum_{a=1}^q\a[a](x)=1\text{ for }\l\text{-a.a. }x\in\mathbb{T}^d\}.
\end{split}
\end{equation}
In what follows we will often write $\a$ instead of $\a\l$. Let eq denote the equidistribution on $\{1,\dots, q\}$. Next we provide the LDP for the KPM.

\begin{prop}\label{LDP}
The measures $\hat\mu_n=\mu_n\circ(\pi_{n})^{-1}$ satisfy a LDP with rate $n^d$ and ratefunction $I-\inf_{\nu\in\PP}I(\nu)$ where
\begin{equation}
I(\nu)=\begin{cases}
  -\b\sum_{a=1}^q\langle J\ast\a[a],\a[a]\rangle+\langle S(\a|\mathrm{eq}),\l\rangle  & \text{if }\nu=\a\l\text{ with }\a\in B\\
  \infty & \text{otherwise. }
\end{cases}
\end{equation}
and the relative entropy is given by $S(\a|\mathrm{eq})=\sum_{a=1}^q\a[a]\log( q\a[a])$.
\end{prop}

Note that we can rewrite the interaction part of the rate function as a 
punishing term for
spatial inhomogeneities and a local term, i.e.
\begin{equation}\label{PunishingLDP}
\begin{split}
I(\nu)=\tfrac{\b}{2}\sum_{a=1}^q&\int du\int dv\bigl[\a[a](u)-\a[a](v)\bigr]^2J(u-v)\cr
&+\int du \bigl[-\b\sum_{a=1}^q\a[a](u)^2+S(\a(u)|\mathrm{eq})\bigr].
\end{split}
\end{equation}
From this we see that global minimizers of $I$ must be flat profiles where $\a[a](u)$ is independent of $u\in\mathbb{T}^d$. 
Indeed, for every $u\in\mathbb{T}^d$
\begin{equation}\label{PunishingLDP2}
\begin{split}
-\b \sum_{a=1}^q\a[a](u)^2+S(\a(u)|\mathrm{eq})\cr
\end{split}
\end{equation}
is the rate function of the mean-field PM given by the Hamiltonian
\begin{equation*}
H_n(\s):=-\frac{1}{n^d}\sum_{x,y\in\D^d_n}1_{\s(x)=\s(y)}, \hspace{1cm}\s\in\O_n
\end{equation*}
and the complete analysis of minimizers is presented in the Ellis-Wang Theorem \cite{ElWa89}. In particular for $q\ge 3$ the model shows a first order phase transition with critical temperature $\b_c(q)=2(q-1)/(q-2)\log(q-1)$. 
The form of its minimizers depends on $\b$ and $q$ but not on $u$ and hence in view of the first summand of \eqref{PunishingLDP}, which punishes spatial inhomogeneities, a global minimizer must be a minimizer of \eqref{PunishingLDP2} equal for every $u\in\mathbb{T}^d$.
%

%

\bigskip
Before we state the main result about GnG of the fuzzy KPM in the next subsection, let us make the following definitions. These are the natural extensions to the Potts situation from the Ising situation in \cite{FeHoMa14}. 
\begin{defn}\label{GnG}
Given any sequence $(\mu_n)_{n\in\N}$ with $\mu_n$ a probability measure on $\O_n$ for every $n\in\N$, define the single-spin conditional probabilities at site 
$u\in\mathbb{T}^d$ as 
\begin{equation}\label{singesitespec}
\g^u_n(\cdot|\a^{(u)}_{n}):=\mu_n\big(\s(\lfloor nu\rfloor)=\cdot\hspace{0.1cm}|\pi^{(u),\s}_n=\a^{(u)}_{n}\big)\hspace{1cm} \a^{(u)}_{n}\in\mathcal{M}^u_n.
\end{equation}
(a) We call a color profile $\a\in B$ \textbf{good} for a sequence of probability measures $(\mu_n )_{n\in\N}$ if there exists a neighborhood $\mathcal{N}_\a\subset B$ of $\a$ 
such that for all $\tilde\a\in\mathcal{N}_\a$ and for all $u\in\mathbb{T}^d$
\begin{equation}
\g^u(\cdot|\tilde\a):=\lim_{n\uparrow\infty}\g^{u}_n(\cdot|\a^{(u)}_{n})
\end{equation}
exists for all sequences $(\a^{(u)}_{n})_{n\in\N}$ with $\a^{(u)}_{n}\in\mathcal{M}^{u}_n$ for every $n\in\N$ such that $\lim_{n\uparrow\infty}\a^{(u)}_{n}=\tilde\a$ in the weak sense. Moreover the limit must be independent of the choice of $(\a^{(u)}_{n})_{n\in\N}$. 

\noindent
(b) A color profile $\a\in B$ is called\textbf{ bad }for $(\mu_n)_{n\in\N}$ if it is not good for $(\mu_n)_{n\in\N}$.

\noindent
(c) $(\mu_n)_{n\in\N}$ is called Gibbs if it has no bad profiles in $B$. 
\end{defn}

\noindent
\textbf{Remarks: }
1) Definition \ref{GnG} (a) implies continuity of $\a\mapsto\g^u(\cdot|\a)$ in the metric $d(\cdot,\cdot)$ defined in \eqref{weakmetric}
 for all $u\in\mathbb{T}^d$ at good profiles. 

\noindent
2) For the KPM $(\mu_n)_{n\in\N}$ all color profiles $\a\in B$ are good since
\begin{equation}\label{SingleSiteCond}
\g^{u}(k|\a)=\frac{\exp(2\b (J\ast\a[k])(u))}{\sum_{l=1}^q\exp(2\b (J\ast\a[l])(u))}
\end{equation}
and hence $(\mu_n)_{n\in\N}$ is Gibbs in the sense of Definition \ref{GnG} (c). 

\noindent
3) Definition \ref{GnG} assigns the notion of Gibbsianness to a sequence of probability measures that live on different spaces. This is different from the notion of Gibbsianness used for example in lattice systems \cite{EnFeSo93,EnFeHoRe02,EnFeHoRe10,Fe05}, but in that respect similar to the definition of Gibbsianness used in the mean-field setting \cite{HaKu04,JaKuRuWe14}. Since there is spatial dependence in our case it makes sense to call the quantity in \eqref{SingleSiteCond} a specification kernel and $\a$ a boundary condition.

\noindent 
4) Definition \ref{GnG} does not consider sequences $(\a^{(u)}_n)_{n\in\N}$ whose weak limit is singular with respect to $\l$. But in Proposition \ref{LDP} we saw that in the thermodynamic limit we can ignore profiles that are singular w.r.t.~the Lebesgue measure or do not lie in the set $B$.

\subsection{The fuzzy Kac Potts model}
Consider the KPM under the local discretisation map $T:\{1,\dots,q\}\mapsto\{1,\dots,s\}$ where $1<s<q$. More precisely, let $R_1,\dots,R_s$ be a partition of $\{1,\dots,q\}$ with $r_i=|R_i|$ and $\sum_{i=1}^sr_i=q$, then $T(a)=i$ if $a\in R_i$. Apply $T$ to all sites simultaneously and consider the fuzzy Kac Potts measure $\mu^T_n:=\mu_n\circ T^{-1}$.
\begin{defn}
We call the generalized fuzzy KPM Gibbs if all profiles $\a\in B$
are good for the sequence $\mu_n^T$.
\end{defn}

In order to determine Gibbsianness of the fuzzy KPM, similar to \eqref{singesitespec}, we 
write for the single-site kernels
\begin{equation}\label{Representation_Kernel0}
\g^{u}_{n,\b,q,(r_1,\dots,r_s)}(k|\nu):=\mu_n^T(\s({\lfloor nu\rfloor})=k|\pi_{n}^{(u),\s}=\nu)
\end{equation}
where $\b$ is the inverse temperature of the KPM and $\nu\in\mathcal{M}^u_n$ with $s$ colors.

\begin{prop}\label{FuzzyKernel}
For each finite $n$ and $u\in\mathbb{T}^d$ we have the representation 
\begin{equation}\label{Representation_Kernel}
\g^{u}_{n,\b,q,(r_1,\dots,r_s)}(k|\nu)= \frac{r_kA^u\big(\b_k(\nu),r_k,\L_k(\nu)\big)}{\sum_{l=1}^sr_lA^u\big(\b_l(\nu),r_l,\L_l(\nu)\big)}
\end{equation}
where $\L_l(\nu)=\{x\in\D_n^d:\text{ }\nu[l]({x/n})=1/n^d\}$, $\b_l(\nu)=\b|\L_l(\nu)|/n^{d}$ and \linebreak
$A^u(\b,r,\L):=\mu_{\L,\b,r}\Bigl(\exp\bigl(2\b(J\ast\pi_{\L}[1])(\frac{\lfloor nu\rfloor}{n})\bigr)\Bigr)$.
Here $\mu_{\L,\b,r}$ denotes the KPM in the subvolume $\L\subset\D^d_n$ with Hamiltonian $$H_{\L}(\s):=-\frac{1}{|\L|}\sum_{x,y\in \L}J(\frac{x-y}{n})1_{\s(x)=\s(y)},$$ inverse temperature $\b$ and $r$ local states.
\end{prop}
%

In view of Proposition \ref{FuzzyKernel} in order to determine GnG of the fuzzy model we must analyse limiting behavior of the constrained KPM $\mu_{\L,\b,r}$ and its continuity properties. The constrained model again satisfies a LDP similar to the one in Proposition \ref{LDP} but now also the spatial structure of the level sets of the conditioning comes into play. We will say that a sequence of diluted sets $\L_n\subset\D_n^d$ converges weakly to the Lebesgue density $\r$ if for all $f\in C(\mathbb{T}^d)$ we have 
$$\frac{1}{n^d}\sum_{x\in \L_n}\d_{x/n}(f)=\frac{1}{n^d}\sum_{x\in \L_n}f(\frac{x}{n})\to\int du\r(u)f(u)$$
as $n\uparrow\infty$ and write $\L_n\Rightarrow\r$.
%
%
\begin{prop}\label{DiLDP} (Diluted version of LDP for empirical color profiles).   
Consider a sequence of diluted sets $\L_n\subset\D_n^d$ with $\L_n\Rightarrow\r$ for some Lebesgue density $\r$ with $N_\r:=\r\l(\mathbb{T}^d)>0$. Denote $\tilde\r(u):=N_\r^{-1}\r(u)$, then the measures $\hat\mu_{\L_n}:=\mu_{\L_n,\b,q}\circ(\pi_{\L_n})^{-1}$ satisfy a LDP with rate $|\L_n|$ and ratefunction $I_{\tilde\r}-\inf_{\nu\in\PP}I_{\tilde\r}(\nu)$ where 
\begin{equation}\label{RateF}
I_{\tilde\r}(\nu)=\begin{cases}
  -\b\sum_{a=1}^q\langle J\ast\tilde\r\a[a],\tilde\r\a[a]\rangle+\langle S(\a|\mathrm{eq}),\tilde\r \l\rangle  & \text{if }\nu[a]=\tilde\r\a[a]\l,\a\in B\\
  \infty & \text{otherwise. }
\end{cases}
\end{equation}
\end{prop}
Note that we can replace the rate $|\L_n|$ by the desired rate $n^d$ since it is arbitrarily close to $|\L_n|N_\r^{-1}$ for large $n$. 
Similar to \eqref{PunishingLDP} we can rewrite $I_{\tilde\r}$ as a sum of two terms, 
i.e.
\begin{equation}\label{DilutedRate}
\begin{split}
I_{\tilde\r}(\nu)=\tfrac{\b}{2}&\sum_{a=1}^q\int du\tilde\r(u)\int dv\tilde\r(v)\bigl[\a[a](u)-\a[a](v)\bigr]^2J(u-v)\cr
&+\int du\tilde\r(u)\Bigl[-b_{\b,\tilde\r,J}(u)\sum_{a=1}^q\a[a]^2(u)+S(\a[\cdot](u)|\mathrm{eq})\Bigr]\cr
\end{split}
\end{equation}
where we defined the site-dependent local temperature as
\begin{equation*}\label{LocalTemp}
\begin{split}
b_{\b,\tilde\r,J}(u):=\b\int dv\tilde\r(v)J(u-v).
\end{split}
\end{equation*}
In this way we have achieved a representation of the large deviation cost 
of profiles of the diluted KPM as an integral over local mean-field PM 
at sites $u$, with $u$-dependent inverse temperatures, and a quadratic 
punishing for spatial inhomogeneity. 
This representation, used for the effective temperatures $\b_l(\nu)$ from Proposition \ref{FuzzyKernel}, will allow us to see that there are no worse conditioning profiles in the fuzzy KPM with class size of at least three than the flat profiles.

%

\bigskip
Let us for the convenience of the reader recall the theorem from \cite{HaKu04} about GnG for the mean-field fuzzy PM which summarizes the precise information on critical parameter values on GnG. Denote by $\b_c(r)$ the inverse critical temperature of the $r$-state mean-field PM.

\begin{thm}\label{HaKu}
Consider the $q$-state mean-field PM at inverse temperature $\b$, 
and let $s$ and $r_1,\dots,r_s$ be positive integers with $1<s<q$ and $\sum^s_{i=1} r_i=q$. Consider the limiting conditional probabilities of the corresponding mean-field fuzzy PM with spin partition $(r_1,\dots,r_s)$.

\medskip

(i) Suppose that $r_i\leq2$ for all $i=1,\dots,s$. Then the limiting conditional probabilities are continuous functions of the empirical mean of the conditioning, for all $\b\geq0$.

\medskip
\noindent  
Assume that $r_i\geq3$ for some $i$ and put $r*:=\min\{r\geq3,r=r_i\text{ for some }i=1,\dots,s\}$, then the following holds. 

\medskip

(ii) The limiting conditional probabilities are continuous for all $\b<\b_c(r*)$.

(iii) The limiting conditional probabilities are discontinuous for all $\b\geq\b_c(r*)$.  
\end{thm}

We now come to the main result, stating that for the fuzzy KPM the critical parameters for GnG are the same as for the mean-field fuzzy PM  if the parameters are such that low temperature Ising classes are avoided. 

\begin{thm}\label{HaKu-Kac}
Consider the $q$-state KPM at inverse temperature $\b$
and let $s$ and $r_1,\dots,r_s$ be positive integers with $1<s<q$ and $\sum^s_{i=1} r_i=q$. 
Consider the limiting conditional probabilities of the corresponding fuzzy KPM with spin partition $(r_1,\dots,r_s)$ where $r*:=\min\{r\geq3,r=r_i\text{ for some }i=1,\dots,s\}$.

\medskip

(i) Suppose that either $\b\leq \b_c(2)$ or that $r_i\not = 2$ for all $i=1,\dots,s$ and $\b<\b_c(r*)$, then the fuzzy KPM is Gibbs. 
The specification kernel is given by 
\begin{equation}\label{Representation_Kernel_Equi}
\lim_{n\uparrow\infty}\g^{u}_{n,\b,q,(r_1,\dots,r_s)}(k|\a_n^{(u)})=\frac{r_k\exp(2\b r_k^{-1}\int dv\r_k(v)J(u-v))}{\sum_{l=1}^sr_l\exp(2\b r_l^{-1}\int dv\r_l(v)J(u-v))}
\end{equation}
when $(\a_n^{(u)})_{n\in\N}$ converges to $\a=(\r_1\l,\dots,\r_s\l)^T$ as defined in Definition \ref{GnG} (a). 

\medskip
(ii) If $r_i\geq3$ for some $i=1,\dots,s$ and $\b\ge\b_c(r_*)$, then the fuzzy KPM is non-Gibbs. 
%
%
%
%
%
%
%
%
%
\end{thm}

\noindent
\textbf{Remarks: }
1) In case $(i)$ the limiting kernels \eqref{Representation_Kernel_Equi} are continuous functions of the conditioning $\a$, as it is explicit from the given expression. 

\noindent
 2) In the mean-field setting, by the fact that for the Ising model phase transitions are of second order, the Ising classes $r_i=2$ can never be a source of discontinuities. This is reflected in part (i) of Theorem \ref{HaKu}. In the fuzzy KPM the situation is potentially 
richer since the Ising classes offer the possibility of a new phenomenon 
related to minimizing profiles which are not spatially homogeneous. This phenomenon, if it occurs, would not be reducible to the mean-field setup. More precisely, for an Ising class, we can re-express the rate function \eqref{DilutedRate} of the diluted LDP
in terms of a $[-1,1]$-valued and site-dependent 
magnetization function $m(u)$ as
\begin{equation*}
\begin{split}
I_{\tilde\r}(m)=\tfrac{\b}{4}&\int du\tilde\r(u)\int dv\tilde\r(v)\bigl[
m(u)-m(v)\bigr]^2J(u-v)+\int du\tilde\r(u)\Phi_u(m(u))\cr
\end{split}
\end{equation*}
where $\Phi_u$ denotes the site-dependent Curie-Weiss Ising rate function obtained 
by substituting the appropriate site-dependent inverse temperature. For certain 
choices of $\tilde \rho$ some sites $u$ can then be made to be in the low-temperature regime 
and others in the high-temperature regime. 
So, for a minimizing magnetization function $m[\tilde \r]$ there is the 
competition between the flatness-imposing term in the double integral and the single-site Curie-Weiss terms which are minimized by $u$-dependent magnetizations. 
This leads us to a non-trivial variational problem and there is a chance for multiple local and global minima in magnetization profile space. In particular there is a possibility for discontinuous behavior of the minimizers under variation of $\tilde\r$ leading to non-Gibbsianness. 
This phenomenon would be a first-order type transition in profile-space, caused genuinely 
by non-homogeneity. 
To decide whether or when it occurs 
clearly deserves more investigation in the future, entering 
the scope of non-homogenous non-convex variational problems.

\noindent
3) Answering a question of a referee, we would like to add the following 
conceptual explanation and provide a small extension.  
For lattice models and models on graphs where a proper DLR formalism is available, the concept of an essential discontinuity of the conditional probabilities of the transformed infinite-volume measure is important, see \cite[Definition 5.13]{Fe05} and \cite{EnErIaKu12}. 
In contrast, for mean-field models or KM, the concept of an essential discontinuity of a {\em limiting specification kernel itself } w.r.t.~the {\em limiting measure $\mu$ itself } is not meaningful.
We must always adopt a sequential view to the approach to the limit to see non-trivial phenomena, as described in Definition \ref{GnG}. 
This is well-established in mean-field models and was successfully adopted for the analysis of a KM in \cite{FeHoMa14}.  

Recall that a $\mu$-essential discontinuity of a function 
is a discontinuity which 
can not be removed by replacing the function by another representative which coincides with it $\mu$-a.s.   

It is not meaningful because the limiting measure  (which in our model a priori is a measure living on the space of Kac-profiles which have densities relative to the Lebesgue measure) tends to be 
a finite combination of Dirac-measures, 
for any inverse temperature.   
This follows in our example from the fact that minimizers of the Kac-rate function 
must be flat and from the Ellis-Wang Theorem for the mean-field empirical color-distribution vectors \cite{ElWa89}. 
More precisely, in our example $\mu$ is supported on a finite set of profiles 
$\a_j$, where $j$ runs from $1$ to at most $q+1$. But for a finitely 
supported measure the a.s.~continuity requirement becomes empty.
Indeed, for any specified values of the limiting kernels 
$\gamma( \cdot |\a_j)$, we could always easily find a continuous interpolating 
function $\a \mapsto \gamma^I( \cdot|\a)$ from the space of profiles to probability vectors, 
which takes the prescribed values. To give such an explicit interpolation, 
take e.g.~the convex combination  $$\gamma^I( \cdot|\a)
:=\sum_j  \Bigl(1 + \sum_{k:k \neq j} \frac{d(\a,\a_j)} {d(\a,\a_k)} \Bigr)^{-1} \gamma( \cdot |\a_j)$$ 
 away from the $\a_j$'s, and $\gamma^I( \cdot|\a_j):=\gamma( \cdot |\a_j)$.

It is however very meaningful here to see whether the finite support of the limiting measure contains a bad configuration in the sense of our Definition \ref{GnG}.
When this is not the case, one would naturally say the model is {\em almost surely Gibbs in the Kac-sense}. 

It is not difficult in our case to conclude 
that this is indeed true throughout
all temperatures, and the fuzzy KPM is almost surely Gibbs in this Kac-sense. 
This is a direct corollary of the present results combined with previous work: 
First one realizes the flatness of profiles in the limiting measure of the KM, from which the corresponding 
profiles of the fuzzy KPM are obtained. But now we are reduced to the mean-field model, for which 
the corresponding result of atypicality of bad configurations was obtained earlier, see \cite{HaKu04}, by computations involving 
the explicit functions appearing as bad configurations, and the typical values
of the mean-field empirical magnetizations,
 following from the Ellis-Wang Theorem.

\section{Proofs}
Let us start with the proofs of the large deviation results. Note that, considering $\L_n\equiv\D_n^d$, Proposition \ref{LDP} is a special case of Proposition \ref{DiLDP}.

\subsection{Proof of Proposition \ref{DiLDP}}
For convenience we write $\mu_{\L_n}$ for $\mu_{\L_n,\b,q}$. 
Let us proceed in two steps. 

\bigskip
\textbf{Step 1: }First we derive the LDP for $J\equiv0$. In this case our Gibbs measure $\mu_{\L_n}$ is just a spatial product measure on $\L_n\subset\D_n^d$ of the equidistribution on $\{1,\dots,q\}$. 
We consider the exponential moment generating function of the color profile at finite discretization $n$ for some $F\in \CC$,
\begin{equation*}
\begin{split}
\mu_{\L_n}[\exp( |\L_n| \pi_{\L_n}(F))]
&=\mu_{\L_n}[\exp( \sum_{a=1}^q\sum_{x\in \L_n}1_{\s(x)=a}F(a,\frac{x}{n}))]\cr
&=\mu_{\L_n}[\prod_{x\in \L_n}\exp( \sum_{a=1}^q1_{\s(x)=a}F(a,\frac{x}{n}))]\cr
&=\prod_{x\in \L_n}\frac{1}{q}\sum_{a=1}^q\exp\big(F(a,\frac{x}{n})\big).
\end{split}
\end{equation*}
Due to spatial independence, we recover the important single-site logarithmic moment generating function $$\G(F(u)):=\log\frac{1}{q}\sum_{a=1}^q\exp(F_a(u)).$$ 
The limit of discretization going to zero for the logarithmic moment generating function of the color profile is given by
\begin{equation*}\label{Limit1}
\begin{split}
\frac{1}{|\L_n|}\log\mu_{\L_n}[\exp(|\L_n|\pi_{\L_n}(F))]=\frac{1}{|\L_n|}\sum_{x\in \L_n}\G(F(\frac{x}{n})) \rightarrow \int du\tilde\r(u)\G(F(u)).\cr 
\end{split}
\end{equation*}
Notice that the diluted rate function 
\begin{equation*}
I_{\tilde\r}(\nu):=
\begin{cases}
  \langle S(\a|\mathrm{eq}),\tilde\r\l\rangle,  & \text{if }\nu=\a\tilde\r\l\text{ with }\a\in B\\
  \infty, & \text{otherwise}
\end{cases}
\end{equation*}
is equivalent to 
\begin{equation*}
\begin{split}
\G_{\tilde\r}^*(\nu):=
\begin{cases}
  \sup_{F\in C}[\nu(F)-\int du\tilde\r(u)\G(F(u))],  & \text{if }\nu=\a\tilde\r\l\text{ with }\a\in B\\
  \infty, & \text{otherwise. }
\end{cases}
\end{split}
\end{equation*}
Indeed, by duality (see also \cite[Lemma 6.2.13]{DeZe10}) it suffices to show that for all $F\in\CC$
\begin{equation}\label{Duality}
\begin{split}
\int du\tilde\r(u)\G(F(u))=\sup_{\nu\in\PP}\Bigl(\nu(F)- I_{\tilde\r}(\nu)\Bigr).
\end{split}
\end{equation}
From this we see that it suffices to take $\nu\in\PP$
with Lebesgue density $\a\tilde\r$ since the r.h.s.~of \eqref{Duality} is equal to minus infinity otherwise. In that case we can write
\begin{equation*}
\begin{split}
\nu(F)- I_{\tilde\r}(\nu)=\int du\tilde\r(u) \Bigl(\langle F(u),\a[\cdot](u)\rangle - S( \a[\cdot](u)|\mathrm{eq}) \Bigr)
\end{split}
\end{equation*}
and the supremum can be considered sitewise. 
Using Jensen's inequality it is easy to see that the supremum is attained in 
$\a[a](u)=\exp{F_a(u)}/\sum_{b=1}^q\exp{F_b(u)}$ and equation \eqref{Duality} is indeed satisfied. That the supremum is achieved follows by convexity (detailed arguments see for example \cite[Lemma 2.6.13]{DeZe10}).
We further note that for continuous $F$ this optimizing profile is even continuous w.r.t.~the spatial variable as well.

\bigskip
\noindent
\textbf{Upper Bound: }
Since 
$\PP$ is compact, all closed sets in $\PP$ are compact and it suffices to consider $K\subset\PP$ compact.
We can assume without loss that $0<\inf_{\nu\in K}I_{\tilde\r}(\nu)$ and hence we can pick $0<a<\inf_{\nu\in K}I_{\tilde\r}(\nu)$. For every $\nu\in K$ there exists a $F_\nu\in C$ such that $\nu(F_\nu)-\int du\tilde\r(u)\G(F_\nu(u))> a$ and the sets 
\begin{equation*}
\begin{split}
U_\nu:=\{\hat\nu\in\PP: \hat\nu(F_\nu)-\int du\tilde\r(u)\G(F_\nu(u))> a\}
\end{split}
\end{equation*}
form an open covering of $K$. Using the Markov inequality we can estimate
\begin{equation*}
\begin{split}
\frac{1}{|\L_n|}&\log\hat\mu_{\L_n}(U_\nu)\cr
&=\frac{1}{|\L_n|}\log\mu_{\L_n}[\exp\bigl({|\L_n|\pi_{\L_n}(F_\nu)}\bigr)>\exp\bigl(|\L_n|( a+\int du\tilde\r(u)\G(F_\nu(u)))\bigr)]\cr
&\leq-a-\int du\tilde\r(u)\G(F_\nu(u))+\frac{1}{|\L_n|}\log\mu_{\L_n}[\exp\bigl({|\L_n|\pi_{\L_n}(F_\nu)}\bigr)]\cr
\end{split}
\end{equation*}
and hence $\limsup_{n\uparrow\infty}\frac{1}{|\L_n|}\log \hat\mu_{\L_n}(U_\nu)\leq -a$ for all $\nu\in K$. Since $K$ is compact it can be covered by a finite number of $U_\nu$ and thus  $\limsup_{n\uparrow\infty}\frac{1}{|\L_n|}\log \hat\mu_{\L_n}(K)\leq -\inf_{\nu\in K}I_{\tilde\r}(\nu)$.

\bigskip
\noindent
\textbf{Lower Bound: }
Let $G\subset\PP$ be open and $G_{\tilde\r\l}$ denote the set of probability measures in $G$ of the form
$\a\tilde\r\l$. If $G_{\tilde\r\l}=\emptyset$, there is nothing to show. Otherwise let $\nu\in G_{\tilde\r\l}$, then  
there exists $\e_1>0$ such that $N_{\e_1}(\nu)\subset G$ and thus using the definition \eqref{weakmetric} we have
\begin{equation*}
\begin{split}
\hat\mu_{\L_n}(G)&\geq\mu_{\L_n}(\pi_{\L_n}\in N_{\e_1}(\nu))=\mu_{\L_n}(d(\pi_{\L_n},\nu)<\e_1)\cr
&=\mu_{\L_n}\Big(\sum_{j=1}^\infty2^{-j}\frac{|(\pi_{\L_n}-\nu)( f_j)|}{1+|(\pi_{\L_n}-\nu)( f_j)|}<\e_1\Big)\cr
&\geq\mu_{\L_n}\Big(\sum_{j=1}^{K(\e_1)}2^{-j}\frac{|(\pi_{\L_n}-\nu)( f_j)|}{1+|(\pi_{\L_n}-\nu)( f_j)|}<\frac{\e_1}{2}\Big)\cr
\end{split}
\end{equation*}
where $K(\e_1)$ is large enough such that $\sum_{j=K(\e_1)+1}^\infty 2^{-j}<\e_1/2$. Further we can estimate 
\begin{equation*}
\begin{split}
\mu_{\L_n}\Big(\sum_{j=1}^{K(\e_1)}2^{-j}\frac{|(\pi_{\L_n}-\nu)( f_j)|}{1+|(\pi_{\L_n}-\nu)( f_j)|}<\frac{\e_1}{2}\Big)&\geq\mu_{\L_n}\Big(\bigcap_{j=1}^{K(\e_1)}\Big\{\frac{|(\pi_{\L_n}-\nu)( f_j)|}{1+|(\pi_{\L_n}-\nu)( f_j)|}<\frac{\e_1}{2}\Big\}\Big)\cr
&=\mu_{\L_n}\Big(\bigcap_{j=1}^{K(\e_1)}\Big\{|(\pi_{\L_n}-\nu)( f_j)|<\e_2\Big\}\Big)\cr
\end{split}
\end{equation*}
where we set $\e_2:=\e_1/(2-\e_1)$. 

In the next step we approximate $\nu$ by probability measures which are flat on a partition of $\mathbb{T}^d$, more precisely, 
we find $\nu^{\rm{flat}}(\nu)$ close to $\nu$ such that $\nu^{\rm{flat}}(\nu)\in M$ where 
\begin{equation*}
\begin{split}
M:=\{\hat\nu\in\PP: \frac{d\hat\nu}{d\l}(u)=\sum_{k=1}^{N'}\hat\a_k\tilde\r(u)1_{C_k}(u)&\text{ for some finite partition } C_k\text{ of }\mathbb{T}^d\cr
&\text{ and some flat colour profile }\hat\a_k \text{ on }C_k\}.
\end{split}
\end{equation*}
Indeed, given any finite partition $(C_k)_{k\in\{1,\dots, N'\}}$ of $\mathbb{T}^d$ where $\tilde\r\l(C_k)>0$ for $k\leq N$ and $\tilde\r\l(C_k)=0$ for $N<k\leq N'$, 
the measure $\nu^{\rm{flat}}(\nu)$ with $d\nu^{\rm{flat}}(\nu)/d\l(u)=\sum_{k=1}^{N'}\a_{k,\nu}\tilde\r(u)1_{C_k}(u)$ where 
\begin{equation*}
\a_{k,\nu}[a]:=
\begin{cases}
  \tilde\r\l(C_k)^{-1}\int_{C_k} du\tilde\r(u)\a[a](u),  & \text{if }\tilde\r\l(C_k)>0\\
  0, & \text{otherwise}
\end{cases}
\end{equation*}
is in $M$.
%
%
Using this, we can approximate for every $j\in\{1,\dots,K(\e_1)\}$
\begin{equation}\label{ColorAndSpace}
\begin{split}
|(\pi^\s_{\L_n}-\nu)( f_j)|&\le |(\pi^\s_{\L_n}-\nu^{\rm{flat}}(\nu))( f_j)|+|(\nu^{\rm{flat}}(\nu)-\nu)( f_j)|
\end{split}
\end{equation}
where for the second summand
\begin{equation*}
\begin{split}
&|(\nu^{\rm{flat}}(\nu)-\nu)( f_j)|\cr
&=\Big|\sum_{a=1}^q\sum_{k=1}^N\big[\a_{k,\nu}[a]\int_{C_k}du\tilde\r(u) f_j(a,u)-\int_{C_k}du\tilde\r(u)\a[a](u) f_j(a,u)\big]\Big|\cr
&=\Big|\sum_{a=1}^q\sum_{k=1}^N\int_{C_k}du\tilde\r(u)\a[a](u)\big[f_j(a,u)-\tilde\r\l(C_k)^{-1}\int_{C_k}dv\tilde\r(v) f_j(a,v)\big]\Big|\cr
&\le\sup_{a\in\{1,\dots,q\}}\sup_{u\in C_k}|f_j(a,u)-\tilde\r\l(C_k)^{-1}\int_{C_k}dv\tilde\r(v) f_j(a,v)|.
\end{split}
\end{equation*}
The $ f_j$ are uniformly continuous and hence it is possible to partition the torus in such a way that for all $a\in\{1,\dots,q\}$ and $j\in\{1,\dots,K(\e)\}$ we have 
\begin{equation}\label{FinePartition}
\begin{split}
\sup_{a\in\{1,\dots,q\}}\sup_{u\in C_k}| f_j(a,u)-\tilde\r\l(C_k)^{-1}\int_{C_k} dv\tilde\r(v) f_j(a,v)|<\frac{\e_2}{3}
\end{split}
\end{equation}
unless $\tilde\r\l(C_k)=0$. Fixing this partitioning, 
for the first summand in \eqref{ColorAndSpace} we have 
\begin{equation*}
\begin{split}
&|(\pi^\s_{\L_n}-\nu^{\rm{flat}}(\nu))( f_j)|\cr
&\leq\sum_{k=1}^{N'} | \sum_{a=1}^q[\frac{1}{|\L_n|}\sum_{x\in \L_n\cap nC_k} f_j(a,\frac{x}{n})1_{\s(x)=a}-\a_{k,\nu}[a]\int_{C_k} du\tilde\r(u)  f_j(a,u)]|\cr
&\leq\sum_{k=1}^{N'}|\sum_{a=1}^q[ \frac{1}{|\L_n|}\sum_{x\in \L_n\cap nC_k}\tilde\r\l(C_k)^{-1}\int_{C_k} dv\tilde\r(v) f_j(a,v)1_{\s(x)=a}-\a_{k,\nu}[a]\int_{C_k} du\tilde\r(u)  f_j(a,u)]|\cr
&\quad+\sum_{k=1}^{N'} | \sum_{a=1}^q \frac{1}{|\L_n|}\sum_{x\in \L_n\cap nC_k} \{f_j(a,\frac{x}{n})-\tilde\r\l(C_k)^{-1}\int_{C_k} dv\tilde\r(v) f_j(a,v)\}1_{\s(x)=a}|\cr
&\le\frac{\e_2}{3}+\Vert f_j\Vert
\sum_{k=1}^{N'}|\sum_{a=1}^q[ \frac{1}{|\L_n|}\sum_{x\in \L_n\cap nC_k}1_{\s(x)=a}-\a_{k,\nu}[a]\tilde\r\l(C_k)]|\cr
\end{split}
\end{equation*}
and thus, 
we can further estimate
\begin{equation*}
\begin{split}
&\mu_{\L_n}\Big(\bigcap_{j=1}^{K(\e_1)}\Big\{|(\pi_{\L_n}-\nu)( f_j)|<\e_2\Big\}\Big)\cr
&\geq\mu_{\L_n}\Big(\bigcap_{k=1}^{N'}\Big\{| \sum_{a=1}^q[ \frac{1}{|\L_n|}\sum_{x\in \L_n\cap nC_k}1_{\s(x)=a}-\a_{k,\nu}[a]\tilde\r\l(C_k)]|<\e_3\Big\}\Big)\cr
&=\prod_{k=1}^{N'} \mu_{\L_n}\Big(\Big\{| \sum_{a=1}^q[ \frac{1}{|\L_n|}\sum_{x\in \L_n\cap nC_k}1_{\s(x)=a}-\a_{k,\nu}[a]\tilde\r\l(C_k)]|<\e_3\Big\}\Big)\cr
\end{split}
\end{equation*}
where $\e_3:=\e_2(3N'\max_{j\in \{1,\dots, K(\e_1)\}}\Vert f_j\Vert)^{-1}$ and we used that $\mu$ is a product measure in the last line. 
Note that for $k\in\{N+1,\dots, N'\}$ the events inside the $\mu_{\L_n}$-measure occur deterministically
for $n$ sufficiently large by the assumption of convergence of the density of the 
set $\L_n$ to zero on those $C_k$ and hence, for large enough $n$, the product restricts to the terms for $k\leq N$. 
For those $k$ 
%
%
%
%
%
%
%
%
introducing the empirical measures $L^\s_{\L_n,k}(a):=|\L_n\cap nC_k|^{-1}\sum_{x\in \L_n\cap nC_k}1_{\s(x)=a}$, then we can further estimate
\begin{equation*}
\begin{split}
&\mu_{\L_n}\Big(\Big\{| \sum_{a=1}^q[ \frac{1}{|\L_n|}\sum_{x\in \L_n\cap nC_k}1_{\s(x)=a}-\a_{k,\nu}[a]\tilde\r\l(C_k)]|<\e_3\Big\}\Big)\cr
&=\mu_{\L_n}\Big(\Big\{\sum_{a=1}^q |\frac{| \L_n\cap nC_k|}{|\L_n|}L^\s_{\L_n,k}(a)-\a_k[a]\tilde\r\l(C_k)|<\e_3\Big\}\Big)\cr
&\geq\mu_{\L_n}\Big(\Big\{|\frac{| \L_n\cap nC_k|}{|\L_n|\tilde\r\l(C_k)}L^\s_{\L_n,k}(a)-\a_k[a]|<\frac{\e_3}{q\tilde\r\l(C_k)},\text{ for all }a\in\{1,\dots,q\}\Big\}\Big).\cr
\end{split}
\end{equation*}
We set $\e_4:=\min_{k\in\{1,\dots,K(\e_1)\}}\e_3/(q\tilde\r \l(C_k))$ and note that $|\L_n\cap nC_k|/(|\L_n|\tilde\r\l(C_k))\to 1$ as $n\uparrow\infty$. Thus we can assume $n$ large enough such that  $\max_{k\in\{1,\dots,K(\e_1)\}}||\L_n\cap nC_k|/(|\L_n|\tilde\r\l(C_k))-1|<\tilde\e<\e_4/2$. Let $\Vert\cdot\Vert_{TV}$ denote the total variational distance of probability measures on $\{1,\dots,q\}$. Then we have
\begin{equation*}
\begin{split}
&\mu_{\L_n}\Big(\Big\{|\frac{|\L_n\cap nC_k|}{|\L_n|\tilde\r\l(C_k)}L^\s_{\L_n,k}(a)-\a_k[a]|<\frac{\e_3}{q\tilde\r\l(C_k)},\text{ for all }a\in\{1,\dots,q\}\Big\}\Big)\cr
&\geq\mu_{\L_n}\Big(\Big\{|\frac{|\L_n\cap nC_k|}{|\L_n|\tilde\r\l(C_k)}L^\s_{\L_n,k}(a)-\a_k[a]|<\e_4,\text{ for all }a\in\{1,\dots,q\}\Big\}\Big)\cr
&\geq\mu_{\L_n}\Big(\Big\{|L^\s_{\L_n,k}(a)-\a_k[a]|<\e_4/2,\text{ for all }a\in\{1,\dots,q\}\Big\}\Big)\cr
&\geq\mu_{\L_n}\Big(\Big\{\Vert L^\s_{\L_n,k}-\a_k\Vert_{TV}<\e_4/4\Big\}\Big).\cr
\end{split}
\end{equation*}
Now we are in the position to apply the lower bound estimate in Sanov's Theorem and write
\begin{equation*}
\begin{split}
\liminf_{n\uparrow\infty}\frac{1}{|\L_n|}\log\hat\mu_{\L_n}(G)&\geq\sum_{k=1}^N \tilde\r\l(C_k)\liminf_{n\uparrow\infty}\frac{1}{|\L_n|\tilde\r\l(C_k)}\log \mu_{\L_n}\Big(\Big\{\Vert L^\s_{\L_n,k}-\a_k\Vert_{TV}<\frac{\e_4}{4}\Big\}\Big)\cr
&\geq-\inf_{\hat\nu\in M_{\e_4}(\nu)}\int du\tilde\r(u)S(\hat\nu(u)|\mathrm{eq})\cr
\end{split}
\end{equation*}
where 
\begin{equation*}
\begin{split}
M_{\e_4}(\nu):=\{\hat\nu\in\PP: \frac{d\hat\nu}{d\l}(u)=\sum_{k=1}^{N'}\hat\a_k\tilde\r(u)1_{C_k}(u)&\text{ for the same partition as } \nu^{\rm{flat}}(\nu)\cr
&\text{ and }\linebreak\max_{k\in\{1,\dots,N'\}}||\hat\a_{k}-\a_k\Vert_{TV}<\frac{\e_4}{4}\}.
\end{split}
\end{equation*}

To finish the proof, we show that
$$\inf_{\nu\in G_{\tilde\r\l}}\inf_{\hat\nu\in M_{\e_4}(\nu)}\int du\tilde\r(u)S(\hat\nu(u)|\mathrm{eq})\le\inf_{\nu\in G_{\tilde\r\l}}\int du\tilde\r(u)S(\nu(u)|\mathrm{eq}).$$
Indeed, since $u\mapsto S(\nu(u)|\mathrm{eq})$ is a convex function, using Jensen's inequality, we have for any $\nu\in G_{\tilde\r\l}$
\begin{equation*}
\begin{split}
\int du\tilde\r(u)S\big(\nu(u)|\mathrm{eq}\big)
&\geq\sum_{k=1}^N\tilde\r\l(C_k)S\Big(\tilde\r\l(C_k)^{-1}\int_{C_k}du\tilde\r(u)\nu(u)|\mathrm{eq}\Big)\cr
&=\sum_{k=1}^N\tilde\r\l(C_k)S(\a_{k,\nu}|\mathrm{eq})=\int du\tilde\r(u)S\big(\nu^{\rm flat}(\nu)(u)|\mathrm{eq}\big)
\end{split}
\end{equation*}
and thus, since $\nu^{\rm flat}(\nu)\in M_{\e_4}(\nu)$, the desired inequality holds.

\bigskip
\textbf{Step 2: }Let us now consider the case with interaction, i.e.~$J\not\equiv 0$. We want to employ Varadhan's Lemma (\cite[Theorem 4.3.1]{DeZe10}) to prove the LDP as in \cite[Theorem 23.19]{Kl08}. 
The conditions in Varadhan's Lemma are indeed satisfied since $J$ is bounded.
%
$\Cox$
\subsection{Proof of Proposition \ref{FuzzyKernel}}
To compute the l.h.s.~of \eqref{Representation_Kernel} write for a fuzzy configuration $\eta\in\{1,\dots,s\}^{\D_n^d\setminus{\lfloor nu\rfloor}}$ where $u\in\mathbb{T}^d$
\begin{equation}\label{CondMeasure}
\begin{split}
\mu^T_{n}(\s({\lfloor nu\rfloor})=k|\s_{\D_n^d\setminus\lfloor nu\rfloor}=\eta)&=\frac{1}{Z_1(\eta)}\sum_{\xi:T(\xi)=(k,\eta)}\mu_{n}(\xi)\cr
&=\frac{1}{Z_2(\eta)}\sum_{\xi:T(\xi)=(k,\eta)}\exp\Big(\b n^d\sum_{a=1}^qF(\pi^\xi_n[a])\Big)
\end{split}
\end{equation}
where $Z_1(\eta)$ and $Z_2(\eta)$ are the appropriate normalization constants. For notational convenience we introduce the notation 
\begin{equation*}
\pi_{n,\L}^\s:=\frac{1}{n^d}\Bigl(\sum_{x\in\L}1_{\s(x)=1}\d_{x/n},\dots,\sum_{x\in\L}1_{\s(x)=q}\d_{x/n}\Bigr)^{T}
\end{equation*} 
for the color profile on $\L\subset\D_n^d$ normalized by $\D_n^d$. In the next step we separate the components in $\pi_n$ corresponding to the site $\lfloor nu\rfloor$. 
We have
\begin{equation*}
\begin{split}
&\sum_{a=1}^qF(\pi^\xi_n[a])=\sum_{a=1}^q\langle J\ast\pi^\xi_n[a],\pi^\xi_n[a]\rangle\cr
&=\sum_{a=1}^q\Big(\langle J\ast\pi_{n,\D_n^d\setminus \lfloor nu\rfloor}^{\xi}[a],\pi_{n,\D_n^d\setminus \lfloor nu\rfloor}^{\xi}[a]\rangle+\tfrac{2}{n^{d}}(J\ast\pi_{n,\D_n^d\setminus \lfloor nu\rfloor}^{\xi}[a])(\tfrac{\lfloor nu\rfloor}{n})1_{\xi(\lfloor nu\rfloor)=a}\Big)\cr
&\hspace{3cm}+n^{-2d}J(0)\cr
&=\sum_{a:T(a)=k}\Big(\langle J\ast\pi_{n,\D_n^d\setminus \lfloor nu\rfloor}^{\xi}[a],\pi_{n,\D_n^d\setminus \lfloor nu\rfloor}^{\xi}[a]\rangle+\tfrac{2}{n^{d}}(J\ast\pi_{n,\D_n^d\setminus \lfloor nu\rfloor}^{\xi}[a])(\tfrac{\lfloor nu\rfloor}{n})1_{\xi(\lfloor nu\rfloor)=a}\Big)\cr
&\hspace{3cm}+\sum_{l\neq k}\sum_{a:T(a)=l}\langle J\ast\pi_{n,\D_n^d\setminus \lfloor nu\rfloor}^{\xi}[a],\pi_{n,\D_n^d\setminus \lfloor nu\rfloor}^{\xi}[a]\rangle+n^{-2d}J(0)
\end{split}
\end{equation*}
where in the last line we used that $T(\xi({\lfloor nu\rfloor}))=k$ assumed in \eqref{CondMeasure}.
The first, third and fourth summand in the last line do not depend on the site $\lfloor nu\rfloor$, in other words, they only depend on the boundary condition $\eta$. Hence in the conditional Gibbs measure \eqref{CondMeasure} corresponding to the above expression the third and fourth summand can be shifted into the normalization constant in the denominator and the remaining two summands can be normalized using the first summand. 
Let us introduce the levelsets of the boundary condition $\L_l(\eta):=\{x\in\D_n^d:\eta(x)=l\}$ then we can write 
\begin{equation*}
\begin{split}
&\Big[\sum_{\xi(\lfloor nu\rfloor):T(\xi(\lfloor nu\rfloor))=k}\sum_{\xi_{\L_k(\eta)}}\exp\Big(\b n^d\sum_{a:T(a)=k}\big(\langle J\ast\pi_{n,\L_k(\eta)}^{\xi}[a],\pi_{n,\L_k(\eta)}^{\xi}[a]\rangle\cr
&\hspace{5cm}+\tfrac{2}{n^{d}}(J\ast\pi_{n,\L_k(\eta)}^{\xi}[a])(\tfrac{\lfloor nu\rfloor}{n})1_{\xi(\lfloor nu\rfloor)=a}\big)\Big)\Big]\cr
&\hspace{4cm}\times\Big[\sum_{\xi_{\L_k(\eta)}}\exp\Big(\b n^d\sum_{a:T(a)=k}\langle J\ast\pi_{n,\L_k(\eta)}^{\xi}[a],\pi_{n,\L_k(\eta)}^{\xi}[a]\rangle\Big)\Big]^{-1}\cr
&=\Big[\sum_{\xi(\lfloor nu\rfloor):T(\xi(\lfloor nu\rfloor))=k}\sum_{\xi_{\L_k(\eta)}}\exp\Big(\sum_{a:T(a)=k}\big((\tfrac{\b|\L_k(\eta)|^2}{n^{d}}\langle J\ast\pi^{\xi}_{\L_k(\eta)}[a],\pi^{\xi}_{\L_k(\eta)}[a]\rangle\cr
&\hspace{6cm}+\tfrac{2\b|\L_k(\eta)|}{n^d}(J\ast\pi^{\xi}_{\L_k(\eta)}[a])(\tfrac{\lfloor nu\rfloor}{n})1_{\xi(\lfloor nu\rfloor)=a})\big)\Big)\Big]\cr
&\hspace{3.5cm}\times\Big[\sum_{\xi_{\L_k(\eta)}}\exp\Big(\tfrac{\b|\L_k(\eta)|^2}{n^{d}}\sum_{a:T(a)=k}\langle J\ast\pi^{\xi}_{\L_k(\eta)}[a],\pi^{\xi}_{\L_k(\eta)}[a]\rangle\Big)\Big]^{-1}\cr
&=\sum_{\xi(\lfloor nu\rfloor):T(\xi(\lfloor nu\rfloor))=k}\mu_{\L_k(\eta),\b\tfrac{|\L_k(\eta)|}{n^d},r_k}\Big[\exp\Big(\tfrac{2\b|\L_k(\eta)|}{n^{d}}(J\ast\pi_{\L_k(\eta)}[\xi(\lfloor nu\rfloor)])(\tfrac{\lfloor nu\rfloor}{n})\Big)\Big]\cr
&=r_k\mu_{\L_k(\eta),\b\frac{|\L_k(\eta)|}{n^d},r_k}\Big[\exp\Big(\tfrac{2\b|\L_k(\eta)|}{n^{d}}(J\ast\pi_{\L_k(\eta)}[1])(\tfrac{\lfloor nu\rfloor}{n})\Big)\Big]\cr
\end{split}
\end{equation*}
as required.
$\Cox$

\subsection{Proof of Theorem \ref{HaKu-Kac}}
\textbf{Part (i): }
First note that a given weakly convergent sequence of boundary conditions $(\nu_n)_{n\in\N}$ in the single-site specification kernel \eqref{Representation_Kernel} is represented in the sequence of level sets $(\L_k(\nu_n))_{n\in\N}$ and in the temperature parameters $(\b_k(\nu_n))_{n\in\N}$ corresponding to the fuzzy classes $k\in\{1,\dots,s\}$. 
For each such fuzzy class $k$ we have a limiting dilution $\r_k$ and limiting inverse temperature $\b N_{\r_k}$   where either $\b\leq \b_c(2)$ or $\b<\b_c(r_*)$ if $r_i\not = 2$ for all $i\in\{1,\dots, s\}$.
In the degenerate case where $\r_k\equiv0$, also $\b N_{\r_k}=0$ and 
$$A^u\big(\b_k(\nu_n),r_k,\L_k(\nu_n)\big)=\mu_{\L_k(\nu_n),\b\frac{|\L_k(\nu_n)|}{n^d},r_k}\Big[\exp\Big(\tfrac{2\b|\L_k(\nu_n)|}{n^{d}}(J\ast\pi_{\L_k(\nu_n)}[1])(\tfrac{\lfloor nu\rfloor}{n})\Big)\Big]$$
converges to $1$ as $n$ tends to infinity as the exponent tends to zero uniformly. If $N_{\r_k}>0$ we can use the LDP given in Proposition \ref{DiLDP}.
We claim that for any such $\r_k$ the rate function \eqref{RateF} is minimized by the flat equidistribution, more precisely the minimizer is given by $\a[\cdot](u)\equiv1/r_k$ away from $\{u\in\mathbb{T}^d: \r_k(u)=0\}$. In order to see this, consider the representation of the rate function given in \eqref{DilutedRate}. Note that, in the second summand, for every $u\in\mathbb{T}^d$
\begin{equation*}
\begin{split}
b_{\b N_{\r_k},\tilde\r_k,J}(u)=\b N_{\r_k}\int dv\tilde\r_{k}(v)J(u-v)=\b\int dv\r_k(v)J(u-v)\le\b
\end{split}
\end{equation*}
which implies, using the Ellis-Wang Theorem \cite{ElWa89} for the mean-field PM and monotonicity of the critical temperatures w.r.t.~the class size, that the equidistribution $\a[\cdot](u)\equiv1/r_k$ is the unique minimizer for every $u$. 
Consequently, since the flat equidistribution also minimizes the first summand in \eqref{DilutedRate}, $(1/r_k)\tilde\r_k\l$ must be the global minimizer of $I_{\tilde\r_k}$.
This implies, that $A^u\big(\b_k(\nu_n),r_k,\L_k(\nu_n)\big)$ 
converges to 
$\exp\bigl(2\b r_k^{-1}\int dv\r_k(v)J(u-v)\bigr)$ as $n$ tends to infinity. Moreover, for any limit profile, the limiting specification kernel of \eqref{Representation_Kernel} is given by \eqref{Representation_Kernel_Equi} and the limit is independent of the approximating sequence. Hence any boundary profile is good according to Definition \ref{GnG} and thus the fuzzy KPM is Gibbs.

\bigskip
\noindent
\textbf{Part (ii): }
First note that at any finite $n$, the single-site conditional probabilities at one given site, depending 
on empirical color profiles away from the single site, are uniquely defined combinatorial objects which 
are given in terms of the elementary formula for conditional probabilities. Hence there is no need and also no freedom 
to talk about different versions of the kernels at finite $n$. 

We show that each bad configuration for the mean-field fuzzy PM provides a bad configuration for the fuzzy KPM when it is interpreted as the spatially homogeneous
(flat) color profile.
In order to prove that a profile $\nu$ is a bad point, according to the Definition \ref{GnG}, it suffices to show that there exist two sequences $\nu_m^+$ and $\nu_m^-$ of conditionings in the fuzzy  KPM which can be realized at some scale $n_m\to\infty$, which 
\begin{enumerate}
\item are both converging to the same limit $\nu$ as $m\to\infty$ weakly, but  
\item which have the property that the limits of the 
kernels $\g^{u}_{n_m,\b,q,(r_1,\dots,r_s)}(k |\nu_m^\pm)$ from formula \eqref{Representation_Kernel0} with the corresponding conditionings $\nu_m^+$ and $\nu_m^-$ exist and are different. 
\end{enumerate}

We will construct those sequences now by a two-step procedure 
as spatial approximants of bad configurations in mean-field.

Bad configurations $\a\in\PP(\{1,\dots,s\})$ for the mean-field fuzzy PM are characterized by the fact that 
for some fuzzy class $r_k$, 
$\b\a[k]=\b_{c}(r_k)$.
Here $\b_c(r_k)$ is the critical temperature parameter where the 
mean-field non-normalized rate function 
\begin{equation}\label{MFRate}
\begin{split}
I^{\rm{MF}}_{\a[k]}(\hat\a):=-\b\a[k]\sum_{a=1}^{r_k}\hat\a[a]^2+S(\hat\a|\mathrm{eq}),\hspace{1cm}\hat\a\in\PP(\{1,\dots,r_k\})
\end{split}
\end{equation}
of the $r_k$-states PM shows a discontinuous (first-order) jump from  
uniqueness to non-uniqueness 
of the global minimizers (for details see \cite{HaKu04}).

Now, consider $\alpha$ such that the set $T \subset \{1,\dots,s\}$ of indices for which $\b\alpha[k]=\beta_c(r_k)$ is non-empty. Let $i$ denote the lowest index in $T$
and pick sequences of length-$s$ probability vectors $\alpha_m^- $
and $\alpha_m^+$ which are given by $\alpha_m^\pm[i]=\alpha[i] \pm1/m$
and $\alpha_m^\pm[l]=\alpha[l] \mp 1/((s-1)m)$ for $l\neq i$ where $m$ tends to infinity.
This construction moves away all conditionings from the critical point. More precisely, for all 
fuzzy classes along the sequences indexed by $m$ the corresponding mean-field model is either 
in the uniqueness regime, $\alpha_m^\pm[k]<\b_c(r_k)$, or in the low-temperature regime, $\alpha_m^\pm[k]>\b_c(r_k)$, for all $k\in\{1,\dots,s\}$ and for all finite sufficiently large $m$.  

The vectors $\alpha_m^\pm$ have to be interpreted as limiting flat profiles $\alpha_m^\pm\l$ in the fuzzy KPM or more precisely as limits of levelsets $\L_n(\alpha_m^\pm[k])\subset\D_n^d\setminus \lfloor nu\rfloor$. For finite $n$, in general, this can only be done approximately. For example we can color $\D_n^d\setminus \lfloor nu\rfloor$ periodically such that every color $k$ appears with frequency $\alpha_m^\pm[k]$ if $\alpha_m^\pm[k]$ is rational. If $\alpha_m^\pm[k]$ is irrational another approximation by rational numbers can be employed. Having done this, we have as $n$ tends to infinity,
$$\frac{1}{n^d}\sum_{x\in \L_n(\alpha_m^\pm[k])}f(\frac{x}{n})\to\alpha_m^\pm[k]\l(f).$$
Now for all $m$, $\alpha_m^-[i]$ is in the uniqueness region of the constrained model and hence, using the diluted LDP as in part (i) of this proof,
\begin{equation}\label{Part2}
\begin{split}
\mu_{\L_n(\alpha_m^-[i]),\b\frac{|\L_n(\alpha_m^-[i])|}{n^d},r_i}\Big[\exp\Big(\tfrac{2\b|\L_n(\alpha_m^-[i])|}{n^{d}}(J\ast\pi_{\L_n(\alpha_m^-[i])}[1])(\tfrac{\lfloor nu\rfloor}{n})\Big)\Big]
\end{split}
\end{equation}
converges to $\exp(2\b r_i^{-1}\alpha_m^-[i])$ as $n$ tends to infinity. This further converges to $\exp(2\b_c(r_i)r_i^{-1})=:\phi^-(r_i)$ as $m$ tends to infinity. On the other hand, for all $m$, $\alpha_m^+[i]$ is in the non-uniqueness region of the constrained model. Since the rate function \eqref{DilutedRate} of the diluted LDP is again given by the mean-field rate function \eqref{MFRate}, the minimizer in the phase-transition regime is given by the Ellis-Wang Theorem, see for example \cite[Theorem 5.3]{HaKu04}. Consequently, \eqref{Part2} where $\alpha_m^-[i]$ replaced by $\alpha_m^+[i]$ converges to 
\begin{equation}\label{Part3}
\begin{split}
\tfrac{1}{r_i}\Big(\exp\bigl(\tfrac{2\tilde\b_m}{r_i}((r_i-1)u(\tilde\b_m,r_i)+1)\bigr)+(r_i-1)\exp\bigl(\tfrac{2\tilde\b_m}{r_i}(1-u(\tilde\b_m,r_i))\bigr)\Big)
\end{split}
\end{equation}
where we abbreviated $\tilde\b_m:=\b\alpha_m^-[i]$ and $u(\b,r)$ is given as the largest solution of the mean-field equation 
\begin{equation*}
\begin{split}
u=(1-\exp(-\b u))/(1+(q-1)\exp(-\b u)),
\end{split}
\end{equation*}
for more details see also \cite{ElWa89}.
For $m$ tending to infinity, using $u(r_i, \b_c(r_i)) = (r_i-2)(r_i-1)^{-1}$, \eqref{Part3} converges to
\begin{equation*}
\begin{split}
\tfrac{1}{r_i}\Big(\exp\bigl(2\b_c(r_i)r_i^{-1}(r_i-1)\bigr)+(r_i-1)\exp\bigl(2\b_c(r_i)r_i^{-1}(r_i-1)^{-1}\bigr)\Big)=:\phi^+(r_i).
\end{split}
\end{equation*}
Let us write $(\alpha_{m,n}^\pm)_{n\in\N}$ for a finite-volume sequence of boundary conditions converging to $\alpha_m^\pm$. Further let $\phi(k)$ denote the limit of \eqref{Part2} where $i$ is replaced by $k\in\{1,\dots,s\}\setminus T$ and note that this limit is independent of the choice of $\pm$. From the previous it follows that there exists a subsequence of volume labels $n_m$ such that for the sequence of profiles 
$\nu_{m}^\pm :=\alpha_{m,n_m}^\pm $ which can be realized at scale $n_m$ we have 
\begin{equation*}\label{Representation_Kernel_Equi1}
\lim_{m\uparrow\infty}\g^{u}_{n_m,\b,q,(r_1,\dots,r_s)}(i|\nu^-_{m})=\frac{r_i\phi^-(r_i)}{\sum_{k\in T\setminus\{i\}}r_k\phi^+(r_k)+\sum_{k\in \{1,\dots,s\}\setminus T}r_k\phi(k)}
\end{equation*}
and 
\begin{equation*}\label{Representation_Kernel_Equi2}
\lim_{m\uparrow\infty}\g^{u}_{n_m,\b,q,(r_1,\dots,r_s)}(i|\nu^+_{m})=\frac{r_i\phi^+(r_i)}{\sum_{k\in T\setminus\{i\}}r_k\phi^-(r_k)+\sum_{k\in \{1,\dots,s\}\setminus T}r_k\phi(k)}.
\end{equation*}
Since $r_i\ge3$ by assumption it is easy to check that $\phi^+(r_i)>\phi^-(r_i)$ and hence 
\begin{equation*}\label{Representation_Kernel_Equi3}
\lim_{m\uparrow\infty}\g^{u}_{n_m,\b,q,(r_1,\dots,r_s)}(i|\nu^-_{m})>\lim_{m\uparrow\infty}\g^{u}_{n_m,\b,q,(r_1,\dots,r_s)}(i|\nu^+_{m}).
\end{equation*}
This concludes the proof.
$\Cox$

\end{document}